            \documentstyle[12pt,graphics]{amsart}
                        \numberwithin{equation}{section}
                        \theoremstyle{plain}
                         \newtheorem{thm}{Theorem}[section]
                          \newtheorem{lem}[thm]{Lemma}
                                \newtheorem{pro}[thm]{Proposition}
                                \newtheorem{cor}[thm]{Corollary}
                        \theoremstyle{definition}
                                \newtheorem{rem}{Remark}[section]
                                        
\newtheorem{defn}{Definition}[section]

\begin{document}

                  \title[Denominators of the Kontsevich integral]
{On denominators of the Kontsevich integral and the universal perturbative
invariant of 3-manifolds}

                \author[ Thang Le ]{Thang T. Q. Le }
\thanks{This work is partially supported by NSF grant DMS-9626404.
This and related preprints can be obtained at {\tt 
http:\linebreak[0]//\linebreak[0]www.\linebreak[0]math.\linebreak
[0]buffalo.\linebreak[0]edu/\linebreak[0]$\sim$letu}
}
\address{Dept. of Mathematics, SUNY at Buffalo, Buffalo, NY 14214, USA }
\email{letu@@newton.math.buffalo.edu}
\begin{abstract} The integrality of the Kontsevich integral
and perturbative invariants is discussed.
We show that the denominator of the  degree $n$ part of the Kontsevich
integral of any knot or link is a divisor of $(2!3!\dots n!)^4(n+1)!$. 
We also show that the denominator of of the degree $n$ part of the universal
perturbative invariant of homology 3-spheres is not divisible by any 
prime greater than $2n+1$.    
                  \end{abstract}
                \maketitle
\addtocounter{section}{-1}                
\newcommand{\G}{\Gamma}
\newcommand{\Q}{{\Bbb Q}}
\newcommand{\bR}{{\Bbb R}}
\newcommand{\CP}{{\cal Q}}\newcommand{\SS}{{\frak S}}
\newcommand{\fS}{\SS}
\newcommand{\cP}{{\cal P}}
\newcommand{\cC}{{\cal B}}
\newcommand{\cH}{{\cal H}}
\newcommand{\CC}{{\cal F}}
\newcommand{\ZCC}{{\cal F_{\Bbb Z}}}
\newcommand{\ZG}{{\cal G_{\Bbb Z}}}
\newcommand{\cA}{{\cal A}}
\newcommand{\hZ}{{\Hat Z}}\newcommand{\cZ}{{\check Z}}
\newcommand{\ZA}{{\cA^{\Bbb Z}}}
\newcommand{\ZP}{{\cP^{\Bbb Z}}}
\newcommand{\R}{{\Bbb R}}\newcommand{\D}{{\cal D}}
\newcommand{\ve}{\varepsilon}
\newcommand{\Z}{{\Bbb Z}}\newcommand{\C}{{\Bbb C}}
\newcommand{\cB}{{\cal B}}\newcommand{\Do}{\overset {\circ}{\cal D}}

        \section{Introduction}
The aim of this paper is to study the denominators of the Kontsevich integral
and the universal perturbative  invariants of homology 3-spheres. 
The Kontsevich integral is a very interesting knot invariant 
(see \cite{Kontsevich}) which
contains in itself all finite type invariants, as well as 
 the Jones polynomial and its numerous generalizations
(see, for example, \cite{BL,BarNatan1}).

The Kontsevich integral has values in a graded algebra of chord diagrams and
 is expressed by a formula involving iterated integrals.
These iterated integrals, when calculated separately, are mostly transcendental
numbers, and can be expressed as rational linear combinations
of multiple zeta values. 
However, if we collect terms together using the relations between chord
diagrams, then 
 J. Murakami and the author showed that the coefficients
 are rational, see \cite{LeMurakami3}. The proof follows  Drinfeld's work
on quasi-Hopf algebras
\cite{Drinfeld2}.

In many problems concerning the Kontsevich integrals and quantum and
perturbative invariants of 3-manifolds, one needs to know the prime factors
of the denominators of the Kontsevich integral.
Examples of such problems
are  the conjecture
7.3 in \cite{LMO} about the relation between perturbative and quantum
 invariants of 3-manifolds and the conjecture of Lawrence in \cite{Lawrence}
about the integrality and  $p$-adic convergence of perturbative invariants
(see also \cite{LawrenceRozansky}).
In this paper we show that the denominator of the
degree $n$ part of the  Kontsevich integral
is a divisor of $(2!\dots n!)^4(n+1)!$, hence it can not have any prime
factor greater than $n+1$. Example shows that the denominators of
the degree 6 part have prime factor 7.
This result plays important role in Ohtsuki's proof  (see \cite{Ohtsuki2})
of the conjecture 7.3 of \cite{LMO} for the $sl_2$ case.

In the joint work \cite{LMO} with  Murakami and Ohtsuki, using the
 Kontsevich integral and the Kirby calculus, we constructed an invariant
$\Omega$ of 3-manifolds  with values in a graded algebra of 3-valent graphs.
We call it a  universal perturbative invariant.
Later
the author in \cite{Le} showed that $\Omega$ is a universal finite type
invariant of integral homology 3-spheres, hence it plays
the role of the Kontsevich integral for homology 3-spheres.
 For the theory of
finite type invariants of homology 3-spheres see \cite{Ohtsuki1,GO}. 
Here we  show that the denominator of the degree $n$ part of $\Omega$
 of rational homology 3-spheres
does not have prime factor greater than $2n+1$.
This result is  closely related to the integrality property
of quantum invariants which says that quantum
invariants (see \cite{Turaev}) of homology 3-spheres at prime
 root of unity are
cyclotomic integer. This had been conjectured by Kontsevich, and was proved in various cases by  Murakami, Takata-Yokota, Masbaum-Wenzl
(see \cite{Murakami,TakataYokota,MasbaumWenzl}).

The idea of the proof of the main result can be explained as follows. 
First we reduce the proof to  establishing the existence of an
{\it associator}, whose denominator has some specific properties.
Associators are solution of a system of equations, important among them are
the so-called pentagon and hexagon equations. The well-known
Knizhnik-Zamolodchikov associator (found by Drinfeld) is not good,
since its coefficients are not even rational. An explicit formula
for this associator is given in \cite{LeMurakami3}. So we search for 
another associator.
Drinfeld used perturbative method to find associators. In this method,
one first finds an associator up to degree $n$, and then  tries to extend
it 
to degree $n+1$. Drinfeld  observed that the obstruction to the extension
is in the cohomology of a certain complex. He then showed that the cohomology
is equal to 0, hence there is no obstruction at all. 
Bar-Natan in \cite{BarNatan2} carried this program over to  the space of chord 
diagrams.
In Drinfeld's and Bar-Natan's papers,
the mentioned  cohomology groups vanish, over the nationals.
We follow Drinfeld's method, trying to solve the hexagons and pentagons
 equations and  keeping track of the denominators.
This leads to the problem of calculating  the cohomology groups over the
 integers.
It turns out
that this cohomology group is a torsion group, annihilated by $(n+1)!$ in
degree $n$. Hence we can estimate the denominator in each step.

There is, however, another difficulty to overcome.
In Drinfeld's and Bar-Natan's papers, in order to solve the hexagon and
pentagon equation, one has to assume some freedom for the so-called $R$-matrix.
But in order to get the Kontsevich integral,  
the $R$-matrix must be fixed and equal to the simplest one.
The purely combinatorial method to solve the hexagon and pentagon equations
does not work if the $R$-matrix is fixed.
We show here that one can still use the perturbative
 method to solve these
equations when  the $R$-matrix is fixed and equal to the simplest one.
This is done by using a result about the uniqueness of the associator up
to gauge transformations (see \cite{LeMurakami3}) and the existence
of a special associator (proved in \cite{Drinfeld2}; the proof
used analysis).

The paper is organized as follows. In \S1 we recall basic definitions
of Chinese character diagrams (chord diagrams). In \S2 we discuss the cobar
 complex of Chinese character diagrams.
Associators are discussed in \S3. We proved the main result about the
existence of an associator with special denominators in \S4. The
results about the  denominators of the Kontsevich integral and $\Omega$
are proved in \S5 and 6. Finally in
\S7 we prove some  technical results.

Acknowledgments: Much of this work was carried out while the author
was visiting the Mathematical Sciences Research Institute in Berkeley.
Research at MSRI is supported in part by NSF grant DMS-9022140.
The author thanks P. Melvin, J. Murakami, T. Ohtsuki,
S. Schack and D. Thurston for helpful discussions.

\section{Chinese character diagrams}

\subsection{Preliminaries}
A {\it uni-trivalent graph} 
is a graph every vertex of which is either univalent
or trivalent. In this paper we always assume that
 each connected component of a uni-trivalent graph contains
at least one univalent vertex. A uni-trivalent graph is
 {\it vertex-oriented} if at each trivalent
 vertex a cyclic order of edges is fixed. Trivalent vertices are also called
 {\em internal vertices}, and  univalent vertices -- {\em external}.
A vertex-oriented uni-trivalent graph is also known a {\em Chinese character}.

Let $X$ be a compact oriented 1-dimensional manifold 
whose components are ordered.
A {\it Chinese character diagram} with support $X$ is the manifold $X$ 
together with a 
vertex-oriented uni-trivalent graph
 whose univalent vertices are on $X$. A Chinese character diagram is
 {\em g-connected}
if the graph is connected.

In all figures the components of $X$ will be depicted by solid lines, 
 the graph  by dashed lines, and
the orientation at every internal vertex is given by
the counterclockwise direction. For this reason the graph is also called
the dashed graph.

Chinese characters and Chinese character diagrams are regarded up to
homeomorphisms preserving  components of the support and orientations
at vertices.

\begin{figure}[htp]
\centerline{
\begin{picture}(60,40)
\put(27.5,17.5){\rotatebox{45}{
\begin{picture}(0,0)
\multiput(0,0)(3,0){7}{\line(1,0){2}}
\end{picture}}}
\put(29.5,17.5){\rotatebox{135}{
\begin{picture}(0,0)
\multiput(0,0)(3,0){7}{\line(1,0){2}}
\end{picture}}}
\multiput(30,18.5)(0,-3){7}{\line(0,1){2}}
\put(30,20){\circle*{2}}
\thicklines
\put(0,0){\vector(1,0){60}}
\end{picture}\hskip 1cm
$=$\hskip 1cm
\begin{picture}(60,40)
\multiput(20,34)(0,-3){12}{\line(0,1){2}}
\multiput(40,34)(0,-3){12}{\line(0,1){2}}
\thicklines
\put(0,0){\vector(1,0){60}}
\end{picture}\hskip 0.6cm
$-$\hskip 0.6cm
\begin{picture}(60,40)
\put(18.5,-2.5){\rotatebox{60}{
\begin{picture}(0,0)
\multiput(0,0)(3,0){15}{\line(1,0){2}}
\end{picture}}}
\put(40,-2.5){\rotatebox{120}{
\begin{picture}(0,0)
\multiput(0,0)(3,0){15}{\line(1,0){2}}
\end{picture}}}
\thicklines
\put(0,0){\vector(1,0){60}}
\end{picture}}
\caption{The STU relation\label{STU}}\end{figure}
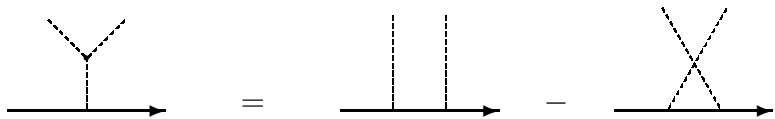

Let $\tilde\cA(X)$ be the vector space over $\Q$ 
spanned by Chinese character diagrams with support $X$.
Let $\cA(X)$ be the quotient space of $\tilde\cA(X)$ by
dividing out  the  STU relation 
shown in Figure \ref{STU}.

The {\it degree} of a Chinese character diagram is  
half the number of vertices of the dashed graph.
Since the relation STU respects the degree, 
there is a grading on $\cA(X)$ induced by this degree. 
We also use  $\cA(X)$ to denote the completion of $\cA(X)$ 
with respect to the degree.

We define a co-multiplication $\hat \Delta$  
in $\cA(X)$  as follows.
A {\it  Chinese character sub-diagram} of a Chinese character diagram $D$ with
 dashed graph $G$
 is any Chinese character
 diagram obtained from $D$ by removing
some (possibly empty) connected components of $G$. The {\it complement
 Chinese character
 sub-diagram} 
of  a Chinese character sub-diagram
$D'$ 
is  the Chinese character  sub-diagram obtained by removing components of $G$ 
which are in $D'$.
We define 
$$\hat \Delta (D)=\sum D'\otimes D''.$$
Here the sum is over all Chinese character sub-diagrams $D'$ of $D$,
 and $D''$ is the complement
of $D'$. This co-multiplication is co-commutative.

\newcommand{\cD}{{\cal D}} 
\newcommand{\g}{{\frak g}}
\newcommand{\cG}{{\cal G}}
\newcommand{\N}{{\Bbb N}}

Suppose that $X,X'$ have  distinguished  components $\ell,\ell'$, and that $X$
 consists of
 loop components only. Let $D\in\cA(X)$ and $D'\in\cA(X')$ be two Chinese
 character diagrams. From each 
 of $\ell,\ell'$ we remove a small arc  which does not contain any
 vertices.
The remaining part of $\ell$ is an arc which we glue  to $\ell'$ in
the place of the removed arc such that the orientations are compatible.
The new Chinese character diagram is called {\it the connected sum of $D,D'$ 
along the 
distinguished components}. It does not depend on the
locations of the removed arcs, which follows from  the STU 
relation and 
the fact that all components of $X$ are loops. The proof is
the same as in case $X=X'=S^1$ as in \cite{BarNatan1}.

In case when $X=X'=S^1$, the connected sum defines a multiplication which
turns  $\cA(S^1)$ into an algebra.

\subsection{Algebra structure} In special cases we can equip $\cA(X)$
with an  algebra structure.
Suppose that  $X$ is $n$ ordered lines on the plane pointing downwards. 
 The space
$\cA(X)$ will be denoted by $\cP_n$; and a connected component of $X$
will be called a {\em string}. 
 If $D_1$ and
$D_2$ are two Chinese character diagrams in $\cP_n$, let  $D_1\times D_2$ be
 the Chinese character diagram obtained  by placing $D_1$ on top of $D_2$.
 The unit
of this algebra is  is the
Chinese character diagram without dashed graph. Let $\cP_0=\Bbb Q$.
It is known that  the algebra $\cP_1$ is commutative (see \cite{BarNatan1}).

The algebra and co-algebra structure
are compatible, and  $\cP_n$ becomes  a Hopf algebra.
It is not hard to see that every primitive element, i.e. element $x$ such that
$\hat \Delta (x)=1\otimes x + x\otimes 1$,  is a linear combination
 of g-connected Chinese character diagrams.

We now introduce a couple of operators acting on $\cP_n$.

 Suppose $D$ is a Chinese character diagram in $\cP_n$ with the
dashed graph $G$.  
Replace the $i$-th string by two strings, the left and  the right, very
close to the old one, and renumber all the strings from left to right.
Attach the graph $G$ to the new set of strings in  the same  way
as in  $D$; this would cause no problem if there is no univalent vertex
on the $i$-th string of $D$. 
If there is a univalent vertex of  $G$  on the $i$-th string then it
yields two possibilities, attaching  to the left or to the right string.
Summing up all $2^m$, where $m$ is the number of univalent vertices
of $D$ on the $i$-th string,  possible Chinese character diagrams of this
 type, we get $\Delta_i(D)\in\cP_{n+1}$. Using linearity we can define
$\Delta_i:\cP_n\to \cP_{n+1}$, for $i=1,2,\dots,n$.

Define $\ve_i$ by $\varepsilon_i(D)=0$ if the Chinese character diagram
 $D\in\cP_n$
 has a univalent  vertex  on the $i$-th string. Otherwise let
 $\varepsilon_i(D)$ be the Chinese character
diagram in $\cP_{n-1}$ obtained  by removing the $i$-th string and renumbering
the remaining strings from left to right.
 We continue $\varepsilon_i$  to a  linear map from $\cP_n$ to $\cP_{n-1}$.

\subsection{Chinese characters}

An {\it $n$-marked Chinese character} $\xi$ is a Chinese character 
whose external vertices are colored by $1,2,\dots,n$.
Two $n$-marked Chinese characters are considered the same if
there is a  homeomorphism between them which  preserves the colors.

Let $\cC_n$ be the vector space over $\Bbb Q$ spanned by all 
$n$-marked Chinese characters subject to the following 
identities:

(1) the antisymmetry identity: $\xi_1+\xi_2=0$, for every two
 Chinese characters $\xi_1$ and $\xi_2$ identical identical everywhere
 except for the orientation at one internal vertex.

(2) the Jacobi identity: $\xi_1=\xi_2+\xi_3$, for every three Chinese
 characters identical outside a ball in which they differ as in Figure
 \ref{IHX}.

\begin{figure}\centerline{
\begin{picture}(60,60)
\multiput(0,20)(3,0){20}{\line(1,0){2}}
\multiput(0,40)(3,0){5}{\line(1,0){2}}
\multiput(15,40)(0,-3){7}{\line(0,-1){2}}
\multiput(30,60)(0,-3){14}{\line(0,-1){2}}
\put(30,3){\makebox(0,0){$\xi_1$}}
\put(15,20){\circle*{2}}
\put(30,20){\circle*{2}}
\end{picture}\hskip 1.5cm
\begin{picture}(60,60)
\multiput(0,20)(3,0){20}{\line(1,0){2}}
\multiput(0,40)(3,0){10}{\line(1,0){2}}
\multiput(30,60)(0,-3){14}{\line(0,-1){2}}
\put(30,3){\makebox(0,0){$\xi_2$}}
\put(30,20){\circle*{2}}
\put(30,40){\circle*{2}}
\end{picture}\hskip 1.5cm
\begin{picture}(60,60)
\multiput(0,20)(3,0){20}{\line(1,0){2}}
\multiput(0,40)(3,0){15}{\line(1,0){2}}
\multiput(45,40)(0,-3){7}{\line(0,-1){2}}
\multiput(30,60)(0,-3){14}{\line(0,-1){2}}
\put(30,3){\makebox(0,0){$\xi_3$}}
\put(30,20){\circle*{2}}
\put(45,20){\circle*{2}}
\end{picture}}
\caption{\label{IHX}}\end{figure}

Now we define a linear mapping $\chi:\cC_n\to\cP_n$ as follows.
Suppose an $n$-marked Chinese character $\xi$ has $k_i$ external vertices
of color $i$.
 There are $k_i!$ ways to put vertices of color $i$ on the $i$-th
string and each of  the $k_1!\dots k_n!$ possibilities gives us a Chinese
character diagram in  $\cP_n$. Summing up all such elements and dividing by
$1/(k_1!k_2!\dots k_n!)$, we get $\chi(\xi)$. It is well-known that
$\chi$  is an
isomorphism between the vector spaces $\cC_n$ and $\cP_n$.
A proof for the case $n=1$ is presented in
\cite{BarNatan1}; the statement itself is Kontsevich's. The proof can be
easily generalized to any $n$. For an explicit description of $\chi^{-1}$,
see \S\ref{chi}.

\section{The cobar complex of Chinese character diagrams}
\subsection{The general complex}

Let us define  $$C^n(\cP)=\cP_n,$$
 $$s_i^n=\ve_i:C^n\to C^{n-1},$$
 and 
$$d^n_i:C^n\to C^{n+1}$$ with $0\le i\le n+1$ by $d^n_i=\Delta_i$
if $1\le i\le n$, and 
$$d^n_0(x)=1\otimes x, \quad \quad d^n_{n+1}=x\otimes 1,$$
where $1\otimes x$ (respectively, $x\otimes 1$) is obtained from $x$ by
 adding a string to the left (respectively, right) of $x$
and renumbering all the strings from left to right.

It was noticed in \cite{BarNatan2} that  $(C^n,d^n_i,s^n_i)$ form  a 
co-simplicial set. It is natural to  consider the following differential
 complex
$$0\to C^0 \overset{d} {\to} C^1 \overset{d} {\to} C^2 \overset{d}
 {\to}\dots C^n \overset{d} {\to} C^{n+1}\dots$$
where 
$$d(x)= \sum_{i=0}^{n+1}(-1)^id^n_i=
1\otimes x -\Delta_1(x)+\Delta_2(x)+\dots+(-1)^n\Delta_n(x)+(-1)^{n+1}
(x\otimes 1).$$
We call it the cobar complex of Chinese character diagrams.
The cohomology of this complex and its subcomplexes will play important
 role.

\subsection{Subcomplexes}
The symmetric group $\SS_n$ acts on $C^n=\cP_n$ on the left by permuting
the strings of the support. An element $x\in\cP_n$ is said to be
 {\it symmetric} if 
$$ x + (-1)^{n(n+1)/2} \sigma (x)=0,$$
where $\sigma$ is the permutation sending 1 to $n$, 2 to $n-1$, 3 to $n-2$,
 etc.

Let $C^n_{sym}(\cP)$ be the subspace of all symmetric elements of $C^n(\cP)$.
 It is easy to see that $(C^*_{sym}, d)$ is a differential subcomplex of
 $(C^*,d)$.

 Let us define the Harrison sub-complex of $\cP_n$ (see  \cite{BarNatan2}).
Let $p,q$ be positive integers. A permutation $\sigma$ of
$\{1,2,\dots, p+q\}$ is called a $(p,q)$-shuffle if $\sigma$ preserves
the order of $\{1,\dots,p\}$ and the order of $\{p+1,\dots,p+q\}$. Define:
$sh_{p,q}\in \Q[\SS_{p+q}]$ by
$$sh_{p,q} = \sum_{\text{all $(p,q)$-shuffles\ } \sigma}(-1)^\sigma \sigma.$$

Let $C^n_{Harr}(\cP)$ be the subset of $C^n(\cP)$ consisting of elements
 $x\in C^n$ such that
$sh_{p,q}(x)=0$ whenever $p+q=n$. As in the usual theory
of Harrison cohomology, one can show that
 $(C^*_{Harr},d)$ is a subcomplex of $(C^n,d)$.

Another way to define the Harrison complex is to use  the Eulerian idempotents.
 It was proved in \cite{Schack} that
there exist idempotents $e_n^{(l)}\in \Q[\SS_n]$, for $l=1,2,\dots,n$ such that
$(e_n^{(l)})^2=e_n^{(l)}$, $e_n^{(l)} e_n^{(k)}=0$ if $k\not = l$, and 
$e_n^{(1)}+\dots e_n^{(n)}=1$. These idempotents, as elements of $\Q[\SS_n]$,
acts on $\cP_n$ as projections. The important point is that all these
idempotents commutes with the differential operator $d$. Hence for each
fixed $l$, we have a subcomplex $(e_n^{(l)} C^n,d)$ of $(C^n,d)$, if we put
$e_n^{(l)}=0$ for  $l>n$.

The subcomplex $(e_n^{(1)}C^n,d)$ is exactly
the above defined Harrison subcomplex.
We record here the formula of $e_n^{(1)}$:
\begin{equation}\label{idempotent}
e_n^{(1)}=\sum_{\sigma\in\SS_n}\frac{(-1)^{r(\sigma)}}
{n {n-1 \choose r(\sigma)}}(-1)^\sigma \sigma,
\end{equation}
where $r(\sigma)$ is the number of $k\in\{1,2,\dots,n-1\}$ such that
$\sigma(k)>\sigma(k+1)$.

\subsection{Non-degeneracy, integral lattices}
An element $x\in\cP_n$ is said to be {\em non-degenerate} if $\ve_i(x)=0$
 for every $i=1,2,\dots,n$. We will be interested in g-connected,
non-degenerate Chinese character diagrams.

Let $\CP^{\Bbb Z}_n$ be the set of all elements in $\cP_n$ which are linear
combinations of g-connected non-degenerate Chinese character
 diagrams with {\em integer} coefficients. Then $\CP^{\Bbb Z}_n$ is a free
$\Z$-module, and let $\CP_n$ be the vector spanned by $\CP^{\Bbb Z}_n$:
$\CP_n=\CP^{\Z}_n\otimes \Bbb Q$.
 From the STU relation one can easily prove the following.
\begin{lem}\label{commutator}
If $x,y$ are  in $\CP^{\Bbb Z}_n$, then so is the commutator $xy-yx$.
\end{lem}

\begin{defn}
Suppose $V$ is a vector space over $\Q$ which contains a fixed integral
lattice $V^{\Z}$, i.e. a $\Z$-module such that $V=V^{\Z}\otimes _{\Z}\Q$.
 We say that an element $x\in V$ has denominator $N$ if
 $Nx$ is in $V^{\Z}$.
\end{defn}

 We will consider $\CP^{\Bbb Z}_n$ the integral lattice of $\CP_n$.

Note that the  operator $d$ preserves $\CP_n$ and $\CP^{\Z}_n$, hence
we can speak about the subcomplexes $(C^*(\CP),d)$ and $(C^*(\CP^{\Z}),d)$,
where $C^n(\CP)=\CP_n$ and $C^n(\CP^{\Z})=\CP^{\Z}_n$.
The latter complex  is a complex over $\Z$. One can also
consider the subcomplexes $(C^*_{sym}(\CP),d)$, $(C^*_{Harr}(\CP),d)$,
$(C^*_{sym}(\CP^{\Z}),d)$ and $(C^*_{Harr}(\CP^{\Z}),d)$ by taking the
intersections with $C^*_{sym}(\cP)$ and $C^*_{Harr}(\cP)$.

It is known that (see \cite{BarNatan2}) the even-dimensional
cohomology groups of  $(C^*_{sym}(\cP),d)$ vanish.
This fact is fundamental in solving the pentagon equation 
(see below) in \cite{BarNatan2,Drinfeld2}.
The proof can be modified easily to show that the
even-dimensional cohomology groups of $(C^*_{sym}(\CP),d)$ vanish.

 Actually, to solve the pentagon equation, one needs the result only for the
 four-dimensional cohomology group. For the purpose of this paper we need
to calculate 
 the four-dimensional cohomology group of 
the  $\Bbb Z$-complex $(C^*_{sym}(\CP^{\Z}),d)$,
 which must have rank 0 as a $\Bbb Z$-module, but may have
some non-trivial torsion part.

Note that $d$ preserves the degree of Chinese character diagrams, hence
the complex $(C^n_{sym}(\CP^{\Z}),d)$ can be decomposed further by 
degree.
The   result, whose  proof will be given in \S\ref{proof},  is

\begin{pro}\label{tech}
The degree $m$ part of the four-dimensional cohomology group of
the complex $(C^n_{sym}(\CP^{\Z}),d)$ is annihilated by 
$2(m+1)!(m!)^2[(m-1)!]^2$.
\end{pro}

\section{ The Drinfeld associator}

\subsection{Associators}
Drinfeld defined associators and $R$-matrix for quasi-triangular quasi-Hopf
algebras, see \cite{Drinfeld1,Drinfeld2}. We recall here the definition,
adapted for the case of Chinese character algebras in 
\cite{LeMurakami2,LeMurakami3} (see also \cite{BarNatan2}).
Note that we don't have any quasi-Hopf algebra here. The algebra
$\cP_n$ will play the role of the $n$-th power of a quasi-Hopf algebra.

Let $r\in\cP_2$ be the Chinese character diagram whose dashed graph is a
line connecting the two strings of the support.
Let $r^{ij}\in\cP_3$ be the Chinese character diagram whose dashed graph is
a line connecting the $i$-th and $j$-th strings of the support.
Define $R=\exp(r/2)\in\cP_2$ and $ R^{ij}=\exp(r^{ij}/2)\in\cP_3$.

\begin{defn}
{\em An associator} is an element $\Phi\in\cP_3$ satisfying
 the following equations: 
 \begin{equation} \Delta_3(\Phi)\times\Delta_1(\Phi)=
(1\otimes \Phi)\times \Delta_2(\Phi)\times(\Phi\otimes 1),\tag{A1}
 \end{equation}
 \begin{equation}
\Delta_1(R)=\Phi^{312}\times
R^{13}\times(\Phi^{132})^{-1}\times R^{23}\times\Phi, \tag{A2}
\end{equation}
\begin{equation}\Phi^{-1}=\Phi^{321},\tag{A3}\end{equation}
\begin{equation}\ve_1(\Phi)=\ve_2(\Phi)=\ve_3(\Phi)=1.\tag{A4}
\end{equation}
\end{defn}
 Here $\Phi^{ijk}$ is the element of
$\cP_3$ obtained from $\Phi$ by  permuting the strings:
the first to the $i$-th, the second to the $j$-th, the third to the
$k$-th. Equation (A1) holds in
$\cP_4$, equations (A2,A3)  in $\cP_3$, and
equation (A4)  in $\cP_2$. There are some redundancy in this system of
 equations. Equation (A1) is known as the pentagon equation, (A2) --
the hexagon equation. Due to (A3), one has that $(\Phi^{132})^{-1}=
\Phi^{231}.$

\begin{rem} Our definition of associator is stricter than that of
\cite{BarNatan2,Drinfeld2} since  the $R$-matrix is $\exp(r/2)$.
\end{rem}

An element in $\cA(X)$ is said to be {\em even} if it is a linear
 combination of 
Chinese character diagrams of even degree. 
It follows from Drinfeld  \cite{Drinfeld2} 
that there is  an {\em even}
associator  $\Phi_{e}$. Evenness has a nice  topological interpretation, see
 \cite{LeMurakami4}. 

\subsection{Symmetric twisting}
Suppose that $F \in\cP_2$ satisfies the following conditions:

\noindent (T1)\quad  $\ve_1(F)=\ve_2(F)=1$,

\noindent (T2)\quad  $F$ is symmetric, i.e. $F$ is a linear combination of 
symmetric Chinese character diagrams.

Then there exists the inverse  $F^{-1}$ in $\cP_2$.
\newcommand{\cF}{{\cal F}}
If $\Phi$ is an element of $\cP_3$, then the element
$$\Phi^F:=[1\otimes
 F]\, \Delta_2(F)\,\,\Phi\,\,\Delta_1(F^{-1})\,[ F^{-1}\otimes
 1]$$
 is said to be obtained from $\Phi$ by {\it twisting via} 
$F$, or by a {\it gauge transformation} (see \cite{Drinfeld1}).

If $\Phi\in\cP_3$ is an associator, i.e. a solution of (A1-A4),
 then it is not difficult to check that $\Phi^F$ is also an associator.
Condition (T2) guarantees that $\Phi^F$ satisfies (A3).

In \cite{LeMurakami3}, following Drinfeld's treatment of the Lie algebra case,
 it was proved that any two associators $\Phi,\Phi'$ are related by a twist
$F$. The proof also shows that if both $\Phi,\Phi'$ are even, then $F$ can be
 chosen to be even. Moreover, if $\Phi,\Phi'$ are ``even associators up to
degree $2k$'', i.e. they satisfy the equations (A1--A4) up to degree $2k$,
then they are related by a twist via an even $F$.

\begin{lem} Suppose that $\Phi\in\cP_3$ is an {\em even} element
which satisfies the hexagon and pentagon equations (A1),
(A2) up to (and including) degree $2k$. Then $\Phi$
also satisfies the same equations up to degree $2k+1$.
\label{even}
\end{lem}
\begin{pf}
The pentagon equation is obviously satisfied, since there is nothing of
odd degree.
For the hexagon equation, we need the fact that there exits an even
associator $\Phi_e$.

 Since both $\Phi$ and $\Phi_e$ are associator up to degree $2k$, there
is an {\em even} element $F\in\cP_2$, satisfying
(T1) and (T2), such that $\Phi^F_e=\Phi$
up to degree $2k$. Since both $\Phi^F_e$ and $\Phi$ don't have
terms of odd degree, they are the same up to degree $2k+1$. The element 
$\Phi^F_e$ is still an associator, hence it satisfies the hexagon equation.
It follows that $\Phi$ satisfies the hexagon equation up to degree $2k+1$.
\end{pf}
\begin{rem} In the proof we used the existence of $\Phi_e$, which was
 established by Drinfeld using analysis. It is still an open problem 
to prove this lemma using only algebra.
\end{rem}

\section{Solving the pentagon and hexagon equations}
\subsection{The existence of a special associator}
To solve the pentagon and hexagon equations we will follow
 \cite{BarNatan2,Drinfeld2}. Lemma \ref{even}  makes the procedure
 much simpler.
Let   $$d_n=(2!3!\dots n!)^4(n+1)!$$

\begin{thm}\label{main}
There exists an {\em even} associator $\Phi\in\cP_3$ of 
 the form $\Phi=\exp(\phi)$, 
where $\phi$ is a linear combination of even, non-degenerate,
g-connected, symmetric Chinese character diagrams. In other words, $\phi$ is
\begin{equation}\label{phi}
\text{even and in} \  C^3_{sym}(\CP).
\end{equation}
Moreover,  the degree $2m$ part of $\phi$ has denominator $d_{2m}$.
\end{thm}
This is the main result. 
Recall that  $\phi\in C^3(\CP)=\CP_3$ is symmetric means that
 $\phi^{321}=-\phi$. The proof of this theorem will occupy the rest of this
 section.

If $\phi$ satisfies (\ref{phi}), then $\Phi=\exp(\phi)$ is even and satisfies
 (A3)
and (A4). There are only the hexagon and pentagon equations to worry about.

For a graded algebra $A$ let $Grad_mA$ (resp. $Grad_{\le m}$) be the
subspace spanned by elements of grading $m$ (of grading $\le m$).

We will solve the pentagon and hexagon equations up to
degree $2m$, and then show that the solution can be extended so that
it solves these equations up to degree $2m+2$.

Suppose that there exists
 $\Phi_{2m}=
\exp(\phi_{ 2m})$
 satisfying (A1) and (A2) up to degree $2m$,
where $\phi_{2m}\in Grad_{\le 2m}\cP_3$  satisfies (\ref{phi})
and the part of degree $2k$ of $\phi_{2m}$ has  denominator $d_{2k}$.
We know that for $m=2$ such a $\Phi_{2m}$ exists, see \cite{BarNatan2}.

 We will find $\phi_{2m+2}=\phi_{2m}+\varphi$, where $\varphi$ is
 of degree $2m+2$ and satisfying (\ref{phi})
 such that $\Phi_{(2m+2)}:=
\exp(\phi_{2m+2})$  satisfies the pentagon and hexagon equations
up to degree $2m+2$. In addition, $\varphi$, and hence $\phi_{2m+2}$,
 has denominator $d_{2m+2}$.

By Lemma \ref{even}, $\Phi_{2m}$ also satisfies the pentagon and hexagon
equations up to degree $2m+1$. Now focus on degree $2m+2$.

\subsection{The hexagon equation}
Let  $\psi$ be the degree $2m+2$ mistake in the hexagon
equation when using $\Phi_{2m}$ in place of $\Phi$, i.e.
\begin{equation}\label{psi}
1+\psi=\Phi_{2m}^{312}\times  (R)^{13}\times (\Phi_{2m}^{-1})^{132}\times
 (R)^{23}\times 
\Phi_{2m} \times  \Delta_1(R^{-1})
\end{equation}
(the equation is taken modulo degree $\ge 2m+3$)

 It is conceivable that with the knowledge of the denominators of the
terms of the right hand side, one can estimate the denominator of $\psi$.
In fact, in \S\ref{proof} we we will show that

\begin{lem}\label{denominator-psi}
The mistake $\psi$  is an element in $\CP_3$ 
 having denominators $[(2m+2)!]^2 d_{2m}$.
\end{lem}

Suppose $\Phi'=\exp(\phi_{2m}+ u)$, where $u$ satisfies (\ref{phi})
 and is of
 degree
$2m+2$. If we replace $\Phi_{2m}$ by $\Phi'$ in Equation (\ref{psi}), then
it's easy to see that the new mistake $\psi'$  is
\begin{equation}\psi'=\psi + u + u^{312} + u^{231}.\label{newpsi}
\end{equation}
Hence if $\psi + u + u^{312} + u^{231}=0$, then $\Phi'$ solves the hexagon
equation up to degree $2m+2$.

The following fact was proved in \cite{BarNatan2,Drinfeld1}.
\begin{lem} The mistake
$\psi$ is totally antisymmetric, i.e. for every permutation $\sigma$ of 3
numbers $\{1,2,3\}$ one has $\sigma(\psi)=sign(\sigma)\psi$.
\end{lem}

Let $u=-\psi/3$, then $\psi + u + u^{312} + u^{231}=0$.
 Hence
 $\bar\Phi:=\exp(\phi_{\le 2m} - \psi/3)$ solves the hexagon equation
up to degree $2m+2$.
 Note that $u$ satisfies (\ref{phi}), is of degree $2m+2$, and has denominator
$3 \times [(2m+2)!]^2 d_{2m}$.

\subsection{The pentagon equation}
Let $\mu$ be the  mistake of degree $2m+2$ in the pentagon equation if
$\Phi$ is replaced by $\bar\Phi$, i.e.
\begin{equation}
1+\mu = \Delta_1(\bar\Phi^{-1})\times \Delta_3(\bar\Phi^{-1})\times (1\otimes
 \bar\Phi)\times 
\Delta_2(\bar\Phi)\times (\bar\Phi
\otimes 1)
\end{equation}
(this equation is taken modulo degrees $\ge 2m+3$).
Again, one can easily estimate the denominator of $\mu$. In \S\ref{proof}
 we will show 

\begin{lem}\label{denominator-mu}
 The mistake  $\mu$ is in $C^4(\CP)=\CP_4$ having denominator  
$3 [(2m+2)!]^2 d_{2m}$.
\end{lem}

If we replace $\bar\Phi$ by $\bar\Phi'=\exp(\phi_{\le 2m}-\psi/3+v)$, where
$v$ satisfies (\ref{phi}) and is of degree $2m+2$, then the new mistake $\mu'$
is
$$
\mu'=\mu -\Delta_1(v)-\Delta_3(v) +\Delta_2(v) + v\otimes 1 + 1\otimes v=
\mu + d v.
$$
The new mistake of the hexagon equation, by (\ref{newpsi}),
 is  $v^{123} +v^{312}+v^{231}$.

Hence, we need to find $v$  of degree $2m+2$ satisfying (\ref{phi}) and the
 following equations
\begin{equation}\label{diff}
dv+\mu=0,
\end{equation} 
\begin{equation}\label{anti}\label{cyclic} 
v^{123} +v^{312}+v^{231}=0.
\end{equation}

The condition (\ref{phi}) says that $v\in C^3(\CP)$ and $v$ is symmetric.
That $v$ is symmetric and satisfying (\ref{cyclic}) is equivalent to
$v$ is symmetric and annihilated by (1,2)- and (2,1)-shuffles. So, we need
to find $v$ in the intersection of $C^3_{sym}(\CP)$ and $C^3_{Harr}(\CP)$
and satisfying $dv= -\mu$.

The following  has been proved in \cite{BarNatan2}.
\begin{lem} The mistake $\mu$ is in the intersection of $C^4_{sym}(\CP)$\and 
$C^4_{Harr}(\CP)$, and $d \mu=0$.
\end{lem}

So $\mu\in C^4_{sym}(\CP)$ and $d\mu=0$.
>From Proposition
\ref{tech} 
 it follows that
there exists  $v'\in C^3_{sym}(\CP)$ having denominator 
$2[(2m+1)!(2m+2)!]^2(2m+3)!$ times that of $\mu$, such that
$dv'=- \mu$.  This $v'$ may not be in $ C^3_{Harr}$, i.e. may not
satisfy (\ref{cyclic}). Put $v= e_3^{(1)}(v')$. Explicitly one has
(see the equation (\ref{idempotent}))
$$v=\frac{2}{3}v' -\frac{1}{3}(v')^{312} -\frac{1}{3}(v')^{231}.$$
Then $v$ is in both $C^3_{sym}(\CP)$ and $C^3_{Harr}(\CP)$.
The commutativity of $e_n^{(1)}$ and $d$ shows  that
 $dv=d(e_3^{(1)}v')= e_4^{(1)} d v' = e_4^{(1)} (-\mu)= -\mu$.
 So $v$ is an element satisfies all (\ref{phi}), (\ref{cyclic})
and (\ref{diff}). Note that
$v$ has denominator $6\times 3 [(2m+1)!]^2[(2m+2)!]^4(2m+3)!d_{2m}$, which is
 a divisor of $d_{2m+2}$ when $m\ge 2$.
 This completes the induction step, and hence the proof
of Theorem \ref{main}.

\section{Denominators of the Kontsevich integral}
\subsection{The Kontsevich integral}
We briefly recall here  the (modification of the) Kontsevich integral for
 framed links and
framed tangles (see \cite{LeMurakami2,LeMurakami3}).
 First we recall the definition of framed tangles.

We  fix an oriented 3-dimensional
 Euclidean space $\bR^3$ with coordinates $(x,y,t)$.
{\it A tangle} is a smooth one-dimensional compact  oriented 
manifold
 $L\subset \bR^3$  lying between two horizontal planes 
$\{t=a\},\{t=b\},
 a<b$ such that all the boundary points are lying on two lines 
$\{t=a,y=0\},
\{t=b,y=0\}$, and at every boundary point $L$ is orthogonal to these 
two planes. These lines are called the top and the bottom lines
 of the tangle.

{\it A normal vector field} on a tangle $L$ is a smooth vector field 
on
$L$ which is nowhere tangent to $L$ (and, in particular, is nowhere
zero) and which is given by the vector $(0,-1,0)$ at every boundary
point. {\it A framed tangle} is a tangle enhanced with a normal vector
field. Two framed tangles are isotopic if they can be deformed by  a
1-parameter family of diffeomorphisms into one another within the 
class
of framed tangles.

Framed oriented links are special framed tangles when there is no 
boundary point.
The empty link, or empty tangle, by definition, is the empty set.

One  assigns a symbol $+$ or $-$ to all the boundary points of a
tangle according to  whether the tangent vector at this point directs
downwards or upwards. Then on the top boundary line of a tangle 
diagram
we have a word  of symbols consisting of $+$ and $-$. Similarly 
on the bottom
boundary line there is another word  of symbols $+$ and $-$.

A {\it non-associative word} on $+,-$ is an element of the  free magma
generated by $+,-$ (see the definition of a free magma in \cite{Serre}).
For a non-associative word one $w$ one defines its length as the
number of its letters. 

 A {\it $q$-tangle} (or 
{\it non-associative}
tangle) is a tangle
together with two non-associative words 
 $w_t,w_b$ such that if we ignore the non-associative structure,
from $w_t,w_b$ we get the words on the top and bottom lines.

If $T_1, T_2$ are framed q-tangles such that $w_b(T_1)=w_t(T_2)$
we can define the product $T=T_1\,T_2$ by placing $T_1$ on 
top of
$T_2$. In this case,
if  $\xi_1\in 
\cA(T_1),
 \xi_2\in \cA(T_2)$ are Chinese character diagrams, then the {\it product} 
$\xi_1\xi_2$
 is a chord diagram in $\cA(T)$ obtained by placing $\xi_1$ on 
top of $\xi_2$.

For any two framed q-tangles $T_1,T_2$ with the same top and bottom lines, 
we can define
their {\it tensor product} $T_1\otimes  T_2$ by putting $T_2$ to 
the right
 of $T_1$. The non-associate structure of the boundaries are the
natural composition of those of $T_1,T_2$.
Similarly, if $\xi_1\in\cP(T_1),\xi_2\in\cP(T_2)$ are chord diagrams, 
then one defines $\xi_1\otimes \xi_2\in \cP(T_1\otimes  T_2)$ 
by the same way.

It is easy to see that every framed  q-tangle $T$
can be obtained from {\em elementary q-tangles},
 using the product and tensor
product. Here an elementary q-tangles is one of the following:

a) a trivial framed q-tangle, i.e. a bunch of vertical lines pointing
downwards, the framing everywhere is given by the vector $(0,-1,0)$.

b) one of the framed q-tangle depicted in Figures \ref{e-tangles}.
Again the framing is given by the vector $(0,-1,0)$.

\begin{figure}
\centerline{
\begin{picture}(240,60)(0,-15)
\put(45,45){\vector(-1,-1){45}}
\put(0,45){\line(1,-1){20}}
\put(25,20){\vector(1,-1){20}}
\put(22.5,-12.5){\makebox(0,0){$X^+$}}
\put(110,45){\line(-1,-1){20}}
\put(65,45){\vector(1,-1){45}}
\put(85,20){\vector(-1,-1){20}}
\put(87.5,-12.5){\makebox(0,0){$X^-$}}
\put(130,45){\vector(0,-1){20}}
\put(175,25){\vector(0,1){20}}
\put(152.5,25){\oval(45,50)[b]}
\put(152.5,-12.5){\makebox(0,0){$Y^+$}}
\put(195,0){\vector(0,1){20}}
\put(240,20){\vector(0,-1){20}}
\put(217.5,20){\oval(45,50)[t]}
\put(217.5,-12.5){\makebox(0,0){$Y^-$}}
\end{picture}}
\caption{}
\label{e-tangles}
\end{figure}

c) the framed q-tangle $T_{w_1,w_2,w_3}$ which has trivial underlying
framed tangle, but the non-associative word of the top is $w_1(w_2w_3)$,
the one of the bottom is $(w_1w_2)w_3$. Here $w_1,w_2,w_3$ are three
arbitrary non-associative words.

d) any of the above with reversing orientation on some of the components.

We will define an invariant $\hZ(T)\in\cA(T)$ for every framed q-tangle
such that
$$
\hZ(T_1T_2)=\hZ(T_1)\hZ(T_2),
$$
$$
\hZ(T_1\otimes T_2)=\hZ(T_1)\otimes \hZ(T_2)
$$
With these requirements, 
One needs only to define  $\hZ$ of the elementary
framed q-tangles. 
For a trivial framed q-tangle $T$ let $\hZ(T)=T$, the Chinese character diagram
in $\cA(T)$ without dashed graph.
Put
$$\hZ(X^\pm)=1 +\frac{\pm 1}{2}  (x)^{(1)}+\dots + \frac{(\pm 1)^n}
{n!2^n}x^{(n)}+\dots,$$
where $x^{(n)}$ is the Chinese character diagrams in $\cA(X^{\pm})$
whose dashed graph consisting of $n$ parallel horizontal
dashed lines connecting
the two solid strings of the support. Let

$$\hZ(T_{w_1w_2w_3})= \Delta^{(|w_1|)}\otimes \Delta^{(|w_2|)}
\otimes \Delta^{(|w_3|)}(\Phi).$$
Here $|w|$ is the length of the word $w$ and $\Delta^{(1)}=id$,
 $\Delta^{(2)}=\Delta$, $\Delta^{(n)}=
\Delta_1\circ\Delta_1\circ \dots \circ\Delta_1$ ($ n-1$ times).
The operation $\Delta^{(n)}$ replaces one string of the support by $n$ strings.
The right hand side of the above equation means that we apply
 $\Delta^{(|w_1|)}$ to the first string, $\Delta^{(|w_2|)}$ to the second,
 and
 $\Delta^{(|w_3|)}$
to the third string of the support of $\Phi$.

If $T'$ is obtained from $T$ by reversing the orientation of
a component $C$, then we put $\hZ(T')=S_C[\hZ(T)]$, where $S_C:
\cA(T)\to\cA(T')$ is the linear map defined as follows.
Suppose  $D\in \cA(T)$ is a Chinese character diagram with $m$ univalent
vertices on $C$. Reversing the orientation of $C$, then multiplying
 by $(-1)^m$,
from $D$ we get $S_C(D)$.

Finally, $\hZ(Y^{\pm})=\sqrt \nu$ where $\nu\in\cP_1$ is obtained from
$S_{C_2}\Phi$ by
identifying the terminating point of the first string with the beginning point
of the second string, and the  terminating point of the second with the
 beginning point
of the third string. Here $C_2$ is the second string. 

These requirements define $\hZ(T)$ uniquely. It is known 
 (see \cite{LeMurakami2,LeMurakami3,BarNatan2}) that
 $\Hat {Z}$ is well-defined and is an isotopy
invariant of framed q-tangles.  
In fact, $\hZ(T)$ is a universal finite type invariant of framed q-tangles.
For more properties of $\hZ$, see \cite{LeMurakami3}.
For a knot $K$, the natural projection from $\cA(S^1)$ to $\cA(S^1)/\approx$
takes $\hZ(K)$ to the Kontsevich integral of $K$. Here $\approx$ is 
the equivalence relation generated by: any Chinese character diagram with
an isolated dashed chord is equivalent to 0. The original Kontsevich integral
is given by an explicit formula (see \cite{Kontsevich}).

It was proved in \cite{LeMurakami3} that $\hZ(T)$ does not depend on the 
associator $\Phi$ if $T$ is a link.

\subsection{Denominators of the Kontsevich integral}

Let $\cA^{\Z}(X)$ be the set of elements in $\cA(X)$ which are linear
combinations of Chinese character diagrams with {\em integer coefficients}.
We consider $\cA^{\Z}(X)$ as the {\em integral lattice} of $\cA(X)$.

\begin{lem} Suppose a q-tangle $T$ is decomposed as $T=T_1 T_2 \dots T_k$, and
that  $x_i\in\cA(T_i), i=1,..,k$ satisfy
\begin{equation}\label{xxx}\text{
the denominator of the degree $m$ part has denominator $d_m$.}
\end{equation}
Then  the element $(x_1\dots x_k)/k!\in\cA(T)$ also satisfies  \ref{xxx}.
In particular,  the product $x_1\dots x_k$ satisfies \ref{xxx}.
\end{lem}
\begin{pf}
An element of degree $n$ in $(x_1\dots x_k)/k!$ has denominator
$k! d_{n_1}\dots d_{n_k}$, where $n_1+\dots+n_k=n$. By Corollary \ref{sao9}
(proved in \S\ref{proof}) $k!d_{n_1}\dots d_{n_k}$ is a divisor of $d_n$.
\end{pf}

By Theorem \ref{main}, $\Phi=\exp(\phi)= \sum_k \phi^k/k!$. Since $\phi$
satisfies (\ref{xxx}), the previous lemma shows that $\Phi$ also
satisfies (\ref{xxx}). The values of $\hZ$ of  elementary q-tangles
satisfy (\ref{xxx}). Hence we have

\begin{thm}
For every framed q-tangle $T$, the degree $m$ part of $\hZ(T)$ has
denominator $d_m=[2!\dots m!]^4(m+1)!$.
In particular, the degree $m$ part of the Kontsevich integral of any knot 
has denominator $d_m$.
\end{thm}

{\em A framed  string link} is a framed tangle containing  no loops
such that the endpoint of any component on the top line projects vertically
to the other endpoint of the same component. In what follows, every
framed string link will be considered as a framed q-tangle, where
the non-associative words of the top and the bottom are  be the same.

Let $L$ be  framed string link with $l$ components
 as follows.
Closing the framed string link in the same way as one closes the braids,
we obtain an $l$ component link. We assume that the framing of the closing
 part is given by the vector $(0,-1,0)$ everywhere. We define the linking
 number of
two components of $L$ to be the linking number of their closures. Similarly,
the self-liking number of a component is the self-linking number of its
 closure.
The linking numbers and self-linking numbers form the linking matrix.

\begin{defn}
An element $x\in\cA(X)$  is said to have i-filter $k$ if it is a linear sum
of Chinese character diagrams, each has at least $k$ internal vertices.
\end{defn}
If $x$ has i-filter $k$, then certainly $x$ has i-filter $k-1$.

\begin{pro}\label{sao7}
Suppose $L$ is a framed string link whose  matrix of linking numbers is 0.
Then $\hZ(L)$ can be represented as a linear sum  $\sum_{k}z_m$,
where $z_m$ has i-filter $m$ and has  denominator not divisible by
any prime greater than $m+2$.
\end{pro}

\begin{pf} It is known that $\hZ(L)\in\cP_l$ is a group-like element
(see \cite{LeMurakami4}). Hence $\hZ(L)=\exp \xi$, where $\xi$ is 
primitive. The element $\xi$ is a linear sum of g-connected 
Chinese character diagrams. Let $\xi= \xi_1 + \xi_2 +\dots$, where
$\xi_k$ has degree $k$. Since the linking matrix is 0, we have $\xi_1=0$.
The g-connectedness implies that $\xi_k$ has i-filter $k-1$.

We have $\xi = \ln \hZ(L)$. From the formula of expansion of  $\ln(1+x)$ and
Lemma \ref{sao2} (proved in \S\ref{proof}),
 it follows that $\xi_k$ has denominator $d_k$.

Return to the formula 
$$\hZ(T)=\exp(\xi_2+\xi_3+\dots).$$ 
Expanding the right hand side, we get a sum of
terms of the form  $x_{n_1}\dots x_{n_k}/k!$. This term has i-filter
$(n_1-1)+ (n_2-1) +\dots (n_k-1)$, and has denominator not divisible
by any prime greater than the maximum $q$ of $\{k,n_1+1,n_2+1,\dots,n_k+1\}$.
Let $m=q-2$. Then $m\le (n_1-1)+ (n_2-1) +\dots (n_k-1)$, since each $n_i$ is
$\ge 2$. Hence this term has i-filter $m$ and denominator not divisible by
any prime greater than $m+2$.
\end{pf}

\section{The universal perturbative invariant of homology 3-spheres}
\subsection{The invariant}
In \cite{LMO} we used the Kontsevich integral and the Kirby calculus to
construct an invariant $\Omega$ of 3-manifolds with values in the
algebra of 3-valent graphs. We  briefly recall here the definition of
$\Omega$.

Let $\D_n$ be the vector space (over $\Q$) generated by
vertex-oriented 3-valent graphs with $2n$ vertices, subject to the 
anti-symmetry and the Jacobi relations.
We don't allow loop components in these 3-valent graphs, so each component
contains at least one, and hence two, vertices. If loops are allowed as
components of vertex oriented 3-valent graphs, we denote the vector
space by $\Do_n$.
For $n=0$ let $\D_n=\Q$.
Set $\D=\prod_{n=0}^\infty \D_n$; the $\D_n$ is the subspace of degree $n$.
Similarly let $\Do=\prod_{n=0}^\infty \Do_n$.
In $\D$ we define the product of two 3-valent graph as their disjoint
union. The unit is the empty graph which is in $\D_0$.

We define first a couple of linear operators.
Suppose  $x\in\cB_l$ is an $l$-marked Chinese character.
 If the number of vertices
of any color is different from $2n$, or if the degree of $x$ is greater
than $(l+1)n$,  put $j'_n(x)=0$. Otherwise,  partitioning
the $2n$ vertices of the same color into $n$ pairs and identifying points of
each pair, from $x$ we get a 3-valent graph which may contains some loops
and which depends on the partition.  Summing
over all possible partitions, we get $j'_n(x)\in\Do$.

There is a map $$pr: \cP_l\to\cA(\sqcup^{l}S^l)$$ defined by identifying
the two endpoints of each (solid) string. 
For $x\in \cA(\sqcup^{l}S^l) $, let $y\in\cP_l$ such that $x=pr(y)$.
 Define 
$j'_n(x)=j'_n(\chi^{-1}(y))$. It is easy to see that $j'_n$ is well-defined.
Note that $j'_n$ lowers the degree by $nl$, and
the values of $j'_n$ is in $Grad_{\le n}\Do$.

Finally, we define $\iota_n(x)$ by replacing every loop in $j'_n(x)$
by $-2n$; the result is in $Grad_{\le n}\D$.

For a framed link $L$ let $\cZ(L)$ be obtained  by 
taking connected sum of $\hZ(L)$ with $\nu$ along every component of $L$. 
It was proved in \cite{LMO} that $\iota_n(\cZ(L)$, where $L$ is a framed link,
does not depend on the orientations of components of $L$, and does not
change under the second Kirby move. Hence the element

$$\Omega_n(L)= \frac{\iota_n(L)}
{\iota_n(U_+)^{\sigma_+}\iota_n(U_-)^{\sigma_-}}\in Grad_{\le n}\D
$$
is an invariant of the 3-manifold $M$ obtained from $S^3$ by surgery along $L$.
Here $U_{\pm}$ are the trivial knots with framing $\pm1$, and $\sigma_\pm$
are the numbers of positive and negative eigenvalues of the linking matrix
 of $L$.

This invariant takes value in $Grad_{\le n}\D$, for a fixed $n$. Without
loss of information, 
 we can combine
all the invariants $\Omega_n$ into one by putting

$$\Omega(M)= 1+\sum_{n=1}^{\infty}Grad_n(\Omega_n(M)).$$

The invariant $\Omega$ has the following important property: for integral
homology 3-spheres, it's degree $n$ part is a universal invariant of
degree $3n$ (see \cite{Le}). For the theory of finite type
invariants of integral homology 3-spheres, see \cite{Ohtsuki1}.
The degrees of finite type invariants of
homology 3-spheres are always multiples of 3 (see \cite{GO}).
 From $\Omega(M)$, together with weight system
coming from Lie algebras (see \cite{BarNatan1}) or from symplectic 
geometry (see \cite{Kapranov,Kontsevich2}) one can construct invariants
 of 3-manifolds
with values in the space of formal power series in one variable.

\subsection{Denominators of $\Omega$}
Note that if $D$ is a Chinese character diagram with less than $2n$ 
external vertices on one solid component, then $j'_n(D)=0$. Since 
$j'_n$ annihilates any Chinese character diagram in $\cA(\sqcup^l(S^1)$ 
of degree $> (l+1)n$, it follows that if $D$ has $>4n$
external vertices on one solid component, then $j'_n(D)=0$.

An element $x$ in $\D$ or in $\Do$ is said to have denominator $N$ if $Nx$ is
a linear combination of 3-valent graphs with {\em integer coefficients}.

\begin{pro}
Suppose that $L$ is a framed link with diagonal linking matrix.
Then $j'_n(\cZ(L))$ has denominator not divisible by any prime greater
 than $2n+1$.
\end{pro}

\begin{pf} Recall that  if $T'$ is obtained from $T$ by increasing the
 framing of a component by 1, then $\hZ(T)$ is obtained from $\hZ(T)$
 by taking connected
sum with $\exp(\theta/2)$ along that component, where $\theta$ is the Chinese
character diagram in $\cA(S^1)$ whose dashed graph is a dashed line
(see \cite{LeMurakami3}).

Suppose that $L$ is obtained from a framed string link $T$ by closing.
In \cite{LeMurakami3} we proved that
$$\hZ(L)= pr [\hZ(T) \times \Delta^{(l)}(\nu)].$$

Changing the framing of each component to 0, from $T$ we get $T'$.
Then $T'$ has 0 linking matrix. One has

$$\hZ(T)= \hZ(T') \times (e^{k_1\theta_1/2}\otimes \dots\otimes
 e^{k_l\theta/2}),$$
where $k_1,\dots,k_l$ are the framings of the components of $T$.
Write  $\hZ(T')= \sum_m z_m$ as in Proposition  \ref{sao7}. Similarly, we
 can write
$\nu=\sum_q y_q$ where $y_q$ has i-filter $q$ and denominator not divisible
by any prime greater than $q+2$. We have 

 $$\hZ(L)= pr [(\sum_m z_m)\times \Delta^{(l)}(\sum_q y_q)
\times (e^k_1\theta_1/2\otimes \dots\otimes e^{k_l\theta/2}].$$

If we expand the right hand side, then we get a sum of terms of the
form 
$$pr[z_m \times  \Delta^{(l)}(y_q)
 \times (\frac{\theta^{n_1}}{n_1! 2^{n_1}}\otimes \dots \otimes  
\frac{\theta^{n_l}}{n_l! 2^{n_l}})].$$
 
Since $j'_n$ annihilates any Chinese character with i-filter $2n+1$, we can 
restrict to the case with $m,q\le 2n$. And since $j'_n$ annihilates
 any Chinese character diagram with with more than $4n$ external vertices
 on one
components, we may assume that all the $n_k$ are less  than or equal to
$2n$. This means the above term has denominator not divisible by
any prime greater than $2n+1$.
It remains to use Lemma \ref{sao5} (proved below) to conclude that $j'_n$ of
this term has denominator not divisible  by any prime greater than $2n+1$.
\end{pf}

\begin{thm}
Suppose that $M$ is a rational homology 3-spheres, i.e. $H_1(M,\Q)=0$.
Then the degree $n$ part of $\Omega(M)$ has denominator
 not divisible by any prime greater than $2n+1$.
\end{thm}

\begin{pf}

It's sufficient to show that $\Omega_n(M)$ has denominator not
divisible by any prime greater than $2n+1$. 
Suppose $M$ is obtained from $S^3$ by surgery along a framed link $L$
with diagonal linking matrix. Then by the previous proposition, 
$j'_n(\cZ(L))$ has denominator not
divisible by any prime greater than $2n+1$ , hence so
do $\iota_n(\cZ(L))$ and $\Omega_n(M)$.

For an arbitrary rational homology 3-sphere $M$, 
 Ohtsuki \cite{Ohtsuki2} showed  there are lens spaces $M_i$ 
of type $(k_i,1)$, $i=1,...,k$ such that
the connected of $M$ and all the $M_i$ is obtained  from $S^3$ by
surgery along a framed link with diagonal linking matrix. 

Since (see \cite{LMO})
 $$\Omega_n(M\#M')=\Omega_n(M) \Omega_n(M'),$$ 
 it's sufficient
to consider the case when $M$ is obtained from $S^3$ by surgery along
a framed link $L$ with diagonal matrix.
\end{pf} 
\subsection{On the map $j'_n$}
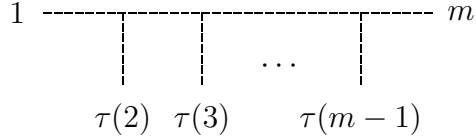
\begin{figure}[htpb]
\centerline{
\begin{picture}(150,45)
\multiput(0,45)(4,0){37}{\line(1,0){3}}
\multiput(30,45)(0,-4){7}{\line(0,-1){3}}
\multiput(60,45)(0,-4){7}{\line(0,-1){3}}
\multiput(120,45)(0,-4){7}{\line(0,-1){3}}
\put(-10,45){\makebox(0,0){1}}
\put(158,45){\makebox(0,0){$m$}}
\put(30,5){\makebox(0,0){$\tau(2)$}}
\put(60,5){\makebox(0,0){$\tau(3)$}}
\put(120,5){\makebox(0,0){$\tau(m-1)$}}
\put(90,25){\makebox(0,0){$\cdots$}}
\end{picture}}
\caption{The definition of $T_{\tau}$}\label{fig.Ttau}
\end{figure}

Let $\cA(m)$, for positive number $m$, be the subspace of $\cC_m$
 spanned by $m$-marked Chinese
characters which have exactly 1 vertex of each color $\{1,2,\dots,m\}$.
Let $\cA(0)=\Q$. 
\label{basaoo}
The symmetric group $\SS_m$ acts on $\cA(m)$ by permuting the colors.
For an element  $\tau$  in the symmetric group $\fS_{m-2}$
acting on the set $\{2,3,\dots,m-1\}$,
let  $T_{\tau} \in \cA(m)$ be the graph shown in Figure \ref{fig.Ttau}.

We define $T_m \in \cA(m)$ by
$$
T_m = \sum_{\tau\in{\fS}_{m-2}}
      \frac{(-1)^{r(\tau)}}{(m-1)\binom{m-2}{r(\tau)}} T_{\tau},
$$
where we denote by $r(\tau)$
the number of $k$ which satisfies $\tau(k)>\tau(k+1)$. The coefficients here
look very similar to those of the Eulerian idempotent $e^{(1)}_{m-1}$, see
(\ref{idempotent}).

There is a shuffle product in  the space $\prod_{m=0}^\infty \cA(m)$
 defined as follows.
Suppose $D$ is a Chinese character in $\cA(m)$, $D'$ -- in $\cA(m')$.
Change the colors of external vertices of $D'$: 1 to $m+1$, 2 to $m+2$,
etc., $m'$ to $m+m'$. The union of $D$ and $D'$ now is an element of
 $\cA(m+m')$. Define
$$D\bullet D' := \sum_{\text{$(m,m')$-shuffles $\sigma \in\SS_{m+m'}$}}
\sigma(D\cup D').$$ 

Let  $T=
\sum_{m=0}^\infty T_m$. Let $T^{\bullet n}$ be the $n$-th power of $T$ in 
the shuffle product. Denote by $T^n_m$ the part of $(T^{\bullet n}/n!)$ 
with $m$ external vertices. In other words,
\begin{equation}
T^n_m= \frac{1}{n!} \sum_{m_1+\dots+m_n=m}T_{m_1}\bullet
 \dots\bullet T_{m_n}.\label{basao4}
\end{equation}
Note that if $m<2n$, then, by definition, $ T^n_m=0$.
The first non-trivial element $T^n_{2n}\in\cA(2n)$ is the following.
 Partition $2n$ points $\{0,1,\dots,2n-1\}$ into $n$ pairs (there are
$(2n-1)!!$ ways to do this), and then connect the two points of each pair
by a dashed line, we get an element of $\cA(2n)$. Summing up all
such possible elements, we get $T^n_{2n}$.

In \cite{LMO} it was proved that the $T^n_m$ satisfy the following important
properties:

1) $T^n_m$ is invariant under the cyclic permutation of
the $m$  external vertices.

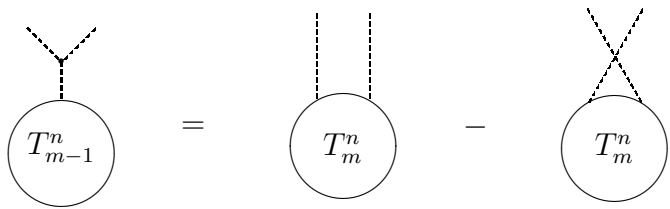
\begin{figure}[htpb]
\centerline{
\begin{picture}(60,80)(0,-38)
\put(28,13){\rotatebox{45}{
\begin{picture}(0,0)
\multiput(0,0)(3,0){6}{\line(1,0){2}}
\end{picture}}}
\put(29,13){\rotatebox{135}{
\begin{picture}(0,0)
\multiput(0,0)(3,0){6}{\line(1,0){2}}
\end{picture}}}
\multiput(30,10)(3,3){3}{\line(1,1){2}}
\multiput(30,15)(0,-3){5}{\line(0,-1){2}}
\put(30,15){\circle*{2}}
\put(30,-20){\circle{40}}
\put(30,-18){\makebox(0,0){$T^n_{m-1}$}}
\put(80,-10){\makebox(0,0){$=$}}
\end{picture} 
\hskip 1.5cm
\begin{picture}(60,80)(0,-38)
\multiput(20,31)(0,-3){11}{\line(0,1){2}}
\multiput(40,31)(0,-3){11}{\line(0,1){2}}
\put(30,-17){\circle{40}}
\put(30,-18){\makebox(0,0){$T^n_m$}}
\put(80,-10){\makebox(0,0){$-$}}
\end{picture}\hskip 1.5cm
\begin{picture}(60,70)(0,-38)
\put(20,-.5){\rotatebox{60}
{\begin{picture}(0,0)
\multiput(0,0)(3,0){14}{\line(1,0){2}}
\end{picture}}}
\put(40,-.5){\rotatebox{120}
{\begin{picture}(0,0)
\multiput(0,0)(3,0){14}{\line(1,0){2}}
\end{picture}}
}
\put(30,-18){\circle{40}}
\put(30,-18){\makebox(0,0){$T^n_m$}}
\end{picture}}
\caption{The dual of the STU relation}\label{dualSTU}
\end{figure}

2) For every $n,m$ we have 
\begin{equation}
 T_m^n - (k\quad\!\! k+1)(T_m^n)  \quad = \quad T_{m-1}^n\star_k Y\tag{*},
\end{equation}
where $(k\quad\!\! k+1)$ is the permutation which interchanges $k$ and
 $k+1$ ($1\le k\le
m-1$),
and $T_{m-1}^n\star_k Y$ denotes the element obtained from $T_{m-1}$
by attaching a Y-shaped graph to the vertex $k$ and then re-numbering
the vertices so that the remaining two vertices of $Y$ are $k$ and $k+1$.
Equation (*) is also presented in 
 Figure \ref{dualSTU}.

For a fixed number $n$, we define a linear map  $j_n:\cA(\sqcup^lS^1)\to
\Do$ 
 as follows. Suppose that $D\in\cA(\sqcup^lS^1)$
is a Chinese character diagram and that  $C$ is a solid circle of $D$
 with $m$ external vertices on it. Number the vertices, beginning
 at any vertex and
following
the orientation of $C$, by 1,2,...,$m$. Now remove the solid circle $C$,
and glue the external vertices to the corresponding vertices of $T_m^n$. 
Do this with all solid circles of
the Chinese character diagram; and  we get $j_n(D)$.

The well-definedness (because of the STU relation)  of this map follows
from equation (*), and this equation can be regarded as the dual to
the STU. 

It follows from the definition that $j_n$ lowers the degree of a
chord diagram by $ln$, where $l$ is the number of solid circles
of $D$.  
 
\begin{lem} Suppose $x\in\cA(\sqcup^l S^1)$ is a Chinese character diagram
of degree less than or equal to $n(l+1)$.
Then $j'_n(x)=j_n(x)$.
\end{lem}
\begin{pf}
Since $\chi$ is an isomorphism, we can assume that $x=pr (\chi( y) )$,
where $y$ is an $l$-marked Chinese character. 
Suppose that $y$ has exactly $2n$ vertices of each color. 
Then $j_n(x)$ is obtained by gluing $T^n_{2n}$ to each group
of vertices of the same  color of $y$, which is, by definition, $j'_n(x)$.
Hence in this case $j_n(x)= j'_n(x)$.

Now suppose one of the color, say $k$, of $y$ has $m$  vertex, where $m\not
 = 2n$.
In this case $j'_n=0$. If $m<2n$, then $j_n(x)=0$, since $T^n_m=0$.
Suppose $m>2n$. Then $T^n_m$ is a linear combination of Chinese character
diagrams, each has at least one  internal vertex, and each  is a
 union of several Chinese character of the form $T_\tau$ (see 
Figure \ref{fig.Ttau}). Hence 
$$\sum_{\sigma\in\SS_m} \sigma (T^n_m)=0,$$
due to the anti-symmetry relation.
On the other hand, the vertices of color $k$ of  the element $\chi(y)$ 
is invariant under actions of $\SS_m$. Hen when we glue
$T^n_m$ to the set of vertices of color $k$, we get 0.
\end{pf}
\begin{pro}Suppose $D$ is a Chinese character diagram in $\cA(\sqcup^l S^1)$.
Then $j'_n(D)$ has denominator not divisible by any prime greater than
$2n+1$.\label{sao5}
\end{pro}

\begin{pf} We can assume that the degree of $D$ is less than or equal to
 $n(l+1)$, since otherwise $j'_n(D)=0$. 
By the previous proposition, $j'_n(D)=j_n(D)$ which is constructed using
$T^n_m$, where $m$ is between $2n$ and $4n$.
 In  (\ref{basao4}), $m_1+\dots+m_n=m$, hence the maximum
of $\{m_1,\dots,m_n\}$ is less than or equal to $2n+2$.
>From the definition it follows that the denominator of these $T_{m_i}$ is not
divisible by any prime greater than $2n+1$.
\end{pf}

\section{Proof of Propositions \protect{\ref{tech}},
\protect{\ref{denominator-psi}} and \protect{\ref{denominator-mu}}}
\label{proof}

\subsection{The cobar complex of Chinese characters} The isomorphism $\chi$
 between $\cP_n$ and $\cB_n$ 
carries the maps $\Delta_i,\ve_i$ over to $\cB_n$. These maps can be
described as follows.

 Suppose $x$ is an $n$-marked Chinese character,
 with $m$ vertices of
color $i$. Here $i\le n$ is a fixed number.
There are $2^m$
ways of partition the set of vertices of color $i$ into
 an ordered pair of subsets, the first and the second subsets.
 For each
such partition, form an $(n+1)$-marked Chinese character by first
changing the color $k$ to $k+1$ for every $k>i$,
then  coloring vertices in the
first subset by $i$, in the second by $i+1$, and leave alone the 
vertices of color $<i$. Summing up, over all possible 
partitions, such $(n+1)$-marked 
Chinese characters, we get $\Delta_i(x)$.

If there is at least one vertex of color $i$, let $\ve_i(x)=0$. Otherwise
$\ve_i(x)$ is the $(n-1)$-marked Chinese character obtained from $x$ by
 changing
the color $k$ to $k-1$ for every $k>i$.

It is easy to check that $\Delta_i,\ve_i$ commute with $\chi$.

For an $n$-marked Chinese character, the number of external vertices is
 called the
{\em  e-grading}, while half the number of all 
vertices is called the {\em degree}.  Note that all the mappings 
$\Delta_i,\ve_i$ preserve  both the e-grading and degree, and $\chi$
 preserves the degree.

The operator $d$ acting on $\cP_n$ is carried 
by $\chi$ over $\cC_n$ to
$d:\cC_n\to\cC_{n+1}$, where
 $$d(\xi)=1\otimes
\xi-\Delta_1(\xi)+\Delta_2(\xi)-\dots+(-1)^n\Delta_n(\xi)+(-1)^{n+1}\xi\otimes
1.$$

Here $1\otimes
\xi$ is the $(n+1)$-marked Chinese character obtained from $x$
 by changing the color $k$ to $k+1$ (for every $k$), and $\xi\otimes 1$
 -- by leaving alone the colors.

\subsection{Subcomplexes}
\newcommand{\CCZ}{\CC^{\Z}}

An element $x\in\cC_n$ is {\em  non-degenerate} if $\varepsilon_i(x)=0$ for
$i=1,2,\dots,n$.
Let $\CC_n$ be the subspace of $\cC_n$ spanned by connected 
non-degenerate Chinese characters.
Let $\CC^{\Z}_{n}$ be the subset of $\CC_n$ consists of elements which
 are linear
combinations of connected non-degenerated Chinese character
 with integer coefficients.
Then 
$$\CC_n=\CCZ_n\otimes_{\Z}\Bbb Q,$$
and we regard $\CC_n$ as the integral lattice of $\CC_n$.
Note that $\chi$ maps $\CC_n$ isomorphically on $\CP_n$. Since $\chi$
 preserves the degree, we can split
this isomorphisms into smaller ones.

Let $\CCZ_n(m)$ (respectively, $\CCZ_n(m,k)$)
 be the subset of $\CCZ_n$ consisting of  elements of degree
$m$ (respectively,  degree $m$ and e-grading $k$). 
Every connected Chinese character of degree $m$ has at most $m+1$
external vertices. Hence
$$\CCZ_n(m)=\oplus_{k=1}^{m+1}\CCZ_n(m,k).$$
 It's important  that the above summation is up to $k=m+1$.
One has 
$$\CC_n(m,k)=\CCZ_n(m,k)\otimes \Bbb Q.$$

Let 
$$C^n(\CCZ(m))=
\CCZ_n(m),$$
 $$C^n(\CCZ(m,k))=\CCZ_n(m,k).$$

\subsection{The inverse of $\chi$}\label{chi}

Let $\cH(m)$ be the  subspace of $\cC_{m+1}$ spanned by $(m+1)$-marked
Chinese characters which have one vertex of each color 1,2,\dots, $m$.
The  other  vertices must have color $m+1$; and we recolor these
ones so that they have color 0.
We define a shuffle product $\cH(m_1)\otimes \cH(m_2)\overset
{\bullet}{\to} \cH(m_1+m_2)$ as follows. Suppose $D,D'$ are
Chinese characters in $\cA(m_1),\cA(m_2)$, respectively.
Change the colors of $D'$: 1 to $m_1+1$, 2 to $m_2+1$, etc, and
$m_2$ to $m_1+m_2$, leaving the 0 color alone. Now $D\cup D'$ is an element
of $\cH(m_1+m_2)$. 
Let 
$$D\bullet D'= \sum_{\text{$(m_1,m_2)$-shuffles $\sigma\in\SS_{m_1+m_2}$}} 
\sigma (D\cup D').$$

\begin{lem}\label{basao8}
If $x\in \prod_{m=0}^\infty \cH(m)$  has denominator 1,
then so does $x^{\bullet n}/n!$.
\end{lem}
\begin{pf} In the shuffle product $x^{\bullet n}$ we have repeated terms 
obtained by
permuting the factor $x$. This cancels the denominator $n!$.
\end{pf}

Consider the element $T_{m+1}\in\cA_{m+1}$ defined in \S\ref{basaoo}.
 Changing
the color $m+1$ to 0, from $T_{m+1}$ we get $t_m\in\cH_m$.
Let 
$$t=t_1+t_2+\dots+t_m+\dots \in \prod_{m=0}^\infty \cH(m).$$
Then $e^t\in\cH$, where we use the shuffle product.
Let $t(m)$ be the part of $e^t$ lying in $\cH_m$ ($m\ge 1$). In  other words,
$$t(m)= \sum_{n=1}^\infty
[\frac{1}{n!} (\sum_{m_1+\dots+m_n=m}t_{m_1}\bullet\dots
\bullet t_{m_n})].$$

The inverse of $\chi$ can be expressed by $e^t$ as follows.
Suppose $D\in\cP_l$ is a Chinese character diagram. Suppose
on the $i$-th string there are $m$ external vertices.
 Remove the $i$-th string,
then glue the vertices of colors 1,2,...,$m$ of  $t(m)$ to
the external  vertices of the $i$-th string, and finally change the color
of the other vertices of $t_m$ from  0 to $i$.
Do this with all the strings. The result is $\chi^{-1}(D)$.
The well-definedness follows from (*). The fact that $\chi(\chi^{-1}(D))=D$
is easy to verify.

\begin{pro}\label{basao}

a) If $x$ is a Chinese character in $\cC_l$ having $m$ external vertices,
then $\chi(x)$ has denominator $m!$.

If $D$ is a Chinese character diagram in $\cP_n$ with $m$ external vertices
then $\chi^{-1}(D)$ has denominator $m!$.

b) If $x$ is a Chinese character in $\CC_l$  having $m$ external vertices,
 where $l\ge 2$,
then $\chi(x)$ has denominator $[(m-1)!]^2.$

If $D$ is a Chinese character diagram in $\CP_l$  having $m$ external vertices,
 where $l\ge 2$,
then $\chi^{-1}(x)$ has denominator $[(m-1)!]^2.$
\end{pro}

\begin{pf} a)
 The first statement follows directly from the definition of 
$\chi$. Let us prove the second statement.
As discussed above, $\chi^{-1}$ can be constricted explicitly using $e^t$.
By Lemma \ref{basao8} $t^n/n!$ has the same denominator as $t$.
>From the formula of  $t_m=T_{m+1}$ in \S\ref{basaoo} it follows that
 has denominator $m!$.

b) follows from the proof of a) with the following observation:
if both $p,q$ are positive integers, then $p!q!$ is a divisor
of $(p+q-1)!(p+q-2)!$.
\end{pf}
\subsection{Some cohomology}
Instead of the complex $(C^*(\CP^{\Z}(m)),d)$, we will study the complex
$C^*(\CCZ(m)),d)$ which is the direct sum of the complexes
$C^*(\CCZ(m,k)),d)$, with $k=1,2,\dots,m+1$:
\begin{equation}\label{decomposition}
C^*(\CCZ(m)),d)= \sum_{k=1}^{m+1}C^*(\CCZ(m,k)),d)
\end{equation}

For a fix number $k$ consider the following $\Z$-complex 
 $E(k)$.
Let $C^n(E(k))$ be the $\Bbb Z$-module spanned by partitions
$\theta_1,\dots,\theta_n$ of $\{1,2,\dots,k\}$, such that each
$\theta_i$ is a non-empty subset of $\{1,2,\dots,k\}$. So if $n>k$ then
$C^n(E(k))=0$. Define $d:C^n(E(k))\to C^{n+1}(E(k))$ by
$$d(\theta_1,\dots,\theta_n)=(d\theta_1,\theta_2,\dots,\theta_n)-
(\theta1,d\theta_2,\theta_3,\dots,\theta_n)+\dots+(-1)^{k-1}(\theta_1,\dots,
d\theta_n),$$
where for a non-empty set $\theta$ we set $d\theta=\sum(\theta',\theta'')$,
the sum is over all possible partition of $\theta$ into an order pair
$\theta',\theta''$ of non-empty subsets.
Actually, $C^*(E(k)),d)$ is the cochain complex of a triangulation of
the $k$-dimensional punctured sphere (see, for example, \cite{LeMurakami3},
\S9.2), and we have
\begin{pro}\label{cube}
One has that $H^k(C^*(E(k)))=\Z$, and $H^l(C^*(E(k)))=0$ for $l\not=k$.
\end{pro}

The symmetric group
$\SS_k$ acts on the left on the complex $C^*(E(k))$ by permuting the
numbers in the partitions. The action is compatible with the operator $d$.
We will see that $C^*(\CCZ(m,k)),d)$ is isomorphic to the tensor
product the complex $C^*(E(k)),d)$ and a $\SS_k$-module which we are going
to describe.

Let $\Gamma^{\Z}(m,k)$ be the subset of 
$\CC_k$  of $\Z$-linear combination of non-degenerate
$k$-marked Chinese character of degree $m$  with exactly one vertex of
each color $\{1,2,\dots,k\}$. And let $\Gamma(m,k)= \Gamma^{\Z}(m,k)
\otimes_{\Z} \Q$.
The symmetric group $\SS_k$ acts on the right  on
$\Gamma^{\Z}(m,k)$  by permuting the colors of external vertices.

\begin{lem} For fixed $m,k$,
there is an isomorphism between complexes 
$(C^*(\CCZ(m,k)),d)$ and
$\Gamma(m,k)\otimes_{\SS_k}(C^*(E(k)),d)$.\label{equiva}
\end{lem}

\begin{pf} Consider a $k$-marked  Chinese character  $\xi$ in  $\G^{\Z}(m,k)$.
The $k$ external vertices are colored by $\{1,2,\dots,k\}$. 
 We
map the element $\xi\otimes(\theta_1,\dots,\theta_n)$ to the element 
$\eta$ of $C^n(\CCZ(m,k))$
obtained from $\xi$ by changing the colors in $\theta_i$ to $i$.
 It can be verified at once that this is an isomorphism
 between the two complexes.
\end{pf}

\begin{lem}\label{order}
The  cohomology groups of the complex
 $\G^{\Z}(m,k)\otimes_{\SS_k}(C^*(E(k)),d)$
is annihilated by $k!$, except for the $k$-th cohomology group.
\end{lem}
\begin{pf}
First let us consider the tensor product over 
$\Bbb Z$: $\G^{\Z}(m,k)\otimes_{\Z}
C^*(E(k)),d)$. This complex has 0 cohomology, except for the cohomology
of dimension $k$, by  Proposition \ref{cube} and the universal coefficient
formula for homology.

There is a natural projection 
$$p_n: \G^{\Z}(m,k)\otimes_{\Bbb Z}C^n(E(k)) \to \G^{\Z}(m,k)
\otimes_{\SS_k}C^n(E(k)).$$
The  kernel $ker(p_n)$ is a finite $\Bbb Z$-module whose cardinality
is a divisor of $k!$.
Consider the short exact sequence of cochain complexes:
$$
0\to ker(p_n)\to \G^{\Z}(m,k)\otimes_{\Bbb Z}C^n(E(k))
\to \G^{\Z}(m,k)
\otimes_{\SS_k}C^n(E(k))
\to 0.$$
By the above observation, the cohomology of the middle complex vanishes
unless in dimension  $k$, while the cohomology of the left one
is annihilated by $k!$. From the long exact sequence derived from
this short exact sequence, we get the lemma.
\end{pf}  

\begin{lem}\label{last}
If $x\in C^4(\CCZ(m,k))$   is symmetric and $dx=0$, then there is
a symmetric  $y\in C^3(\CC(m,k))$,
 such that $dy=x$. Moreover, $y$ has denominator $2\, k!$.
\end{lem}
\begin{pf}
By Lemma \ref{equiva} we may suppose that $x\in \G^{\Z}(m,k)
\otimes_{\SS_k}C^4(E(k))$. Consider two cases: $k=4$ and $k\not=4$.

Suppose $k=4$. Then $x$ is equal to $\gamma\otimes e$, where
$\gamma\in\G^{\Z}(m,4)$ and $e=(\{1\},\{2\},\{3\},\{4\})$, a partition
of $\{1,2,3,4\}$.

Let $$y'= (12,3,4)-(2,13,4)+(2,3,14)-(2,34,1)+(24,3,1)-(4,23,1),$$
where, for example, $(2,13,4)$ means the element $\gamma\otimes
(\{2\},\{1,3\},\{4\})$ in $\G^{\Z}(m,4)\otimes_{\SS_4}C^3(E(4))$.
One can readily check that $dy'=x - x^{4321} = 2x$.

Let $y=(1/4) [y'-(y')^{321}]$; we see that $y$ is symmetric and $dy=x$.

Now consider the case $k\not= 4$. By Lemma \ref{order} there is
$y'\in\G^{\Z}(m,k)\otimes_{\SS_k}C^3(E(k))$ such that $dy'=(k!)x$.
Now $y=(1/2\, k!)[y'-(y')^{321}]$ is the element to find. 
\end{pf}

\subsection{Proof of Proposition \protect{\ref{tech}}}
Suppose $\xi\in C^4_{sym}(\CP^{\Z})$ is of
degree $m$ and $d\xi=0$. Then $\xi'=\chi^{-1}(\xi)$ is symmetric and in
$C^4(\CCZ(m))$. In addition, $\chi(\xi)$ has denominator $(m!)^2$, since $\xi$,
as an element of degree $m$, can have at most $m+1$ external vertices.
(Here we use the fact that $\xi$ is a linear combination of g-connected
elements, and Proposition \ref{basao}, part b)).

Using the decomposition (\ref{decomposition}), we can assume that $\xi'$ is in
$C^4(\CC(m,k))$ with some $k\le m+1$. Besides $\xi'$ is symmetric,
 has denominator $(m!)^2$ and $d\xi'=0$. From Lemma \ref{last} it follows that
there is $y\in C^3(\CC(m,k))$, symmetric, with denominator $2(m!)^2k!$
such that $dy=\xi'$. Since $k\le m+1$, the number $2m!k!$
is a divisor of $2m!(m+1)!$.
Now let $\eta=\chi^{-1}y$. Then $d\eta=\chi^{-1}(\xi')=\xi$, and $\eta$ is
symmetric. By Proposition \ref{basao}, $\eta$ has denominator $2(m!)^2(m+1)!$.
This completes the proof of Proposition \ref{tech}. 

\subsection{On the numbers $d_n$} We will write  $p \lhd q$ if $p$ is
a divisor of $q$.

\begin{lem}\label{sao2}

a) For any positive integers $p,q$ one has \quad $(p+q)!\, d_p d_q \lhd\,
 d_{p+q}$.

b) For any positive integers $p,q$ one has \quad $ d_p d_q \lhd\,  2d_{p+q-1}$.

c) For any integers $p,q \ge 2$ one has \quad $ d_p d_q \lhd\, 96 d_{p+q-2}$.
 
\end{lem}

\begin{pf}
We will prove a). The others can be proved in a similar way.

Note that $d_{p+1}= d_p \times [(p+1)!]^3 (p+2)!$. We use induction on $q$.
The statement is true for $q=1$. Suppose it has been true for $q$.
Then, 
$$(p+q+1)!\, d_p d_{q+1}= (p+q+1)! d_p d_q [(q+1)!]^3 (q+2)!
$$ 
By the induction hypothesis, the latter is a divisor of
$(p+q+1) d_{p+q} [(q+1)!]^3 (q+2)!$

Since $(p+q+1)(q+1)! \lhd (p+q+1)!$,  the number 
$(p+q+1) d_{p+q} [(q+1)!]^3 (q+2)!$
 is a
divisor of $d_{p+q}[(p+q+1)!]^3 (p+q+2)!$ which is $d_{p+q+1}$. This completes
the proof of a).
\end{pf} 

\begin{cor} \label{sao9}
 For every
positive integers $n_1,n_2,\dots n_k$, the number  
 $ k! d_{n_1} d_{n_2} \dots d_{n_k}$ is a divisor of $ d_{n_1+n_2+\dots n_k}$.

\end{cor}
\begin{pf}
>From part a) of the previous lemma, one has that 
$ (n_1+n_2+\dots + n_k)! d_{n_1} d_{n_2} \dots d_{n_k}
 \lhd d_{n_1+n_2+\dots n_k}$.
It remains to notice that $ k \le (n_1+n_2+\dots +n_k)$, since
each $n_i$ is a positive integer.
\end{pf}

We recall here the  Campbell-Hausdorff formula. Let $B={\Bbb Q}\langle
\langle x^{(1)}, x^{(2)},\dots ,x^{(l)}\rangle\rangle$ be the algebra of
 formal power series in $l$
non-commuting variables.
This algebra is graded by the degree of monomials in $x^{(j)}, j=1,2...,l$. 
The free Lie 
algebra $\cal L$ over $\Bbb Z$ (the set of integers) generated by $x^{(j)},
 j=1,2...,l$ is a subset of $B$, and an element in $B$
 is called a {\it  Lie polynomial} if it is
a $\Bbb Q$-linear combination of elements in $\cal L$. 

The  Campbell-Hausdorff formula says that
$$\exp(x^{(1)})\dots \exp(x^{(l)}) =\exp[\sum_{k=1}^\infty f_k(x^{(1)},\dots,
x^{(l)})],$$
where $f_k$ is a homogeneous Lie polynomial of total degree $k$.
Moreover from the Dynkin's form of the Campbell-Hausdorff formula
 (see, for example \cite{Serre}), one can  easily see  that each $f_k$
is a $\Bbb Q$-linear combination of elements in $\cal L$ whose
coefficients have denominators  $(k!)^2$.

\subsection{Proof of Lemmas \ref{denominator-psi} and \ref{denominator-mu}}
We have that (modulo part of degree $\ge 2m +3$)
$$
1+\psi =\exp(\phi_{2m}^{312})\times \exp(\frac{r^{13}}{2})
\times \exp(-\phi_{2m}^{132})\times \exp(\frac{r^{23}}{2})
\times \exp(\phi_{2m})\times \exp(\frac{-r^{13}-r^{23}}{2}).
$$ 

 Hence
\begin{equation}\label{lasst}
\psi=\sum_k f_k(\phi_{2m}^{312},\frac{r^{13}}{2},-\phi_{2m}^{132},
\frac{r^{23}}{2},
 \phi_{2m},\frac{-r^{13}-r^{23}}{2}).
\end{equation}
Recall that $\phi_{2m}$ is a sum of even degree parts, the part of
degree $2k$ has denominator $d_{2k}$. The element
$r^{ij}/2$ has degree 1, and has denominator $d_1=2$.
 Replacing $\phi_{2m}$ by the sum of its degree $2k$ parts, $k=1,\dots,m$,
then expanding the right hand side, we see that $\psi$ is a sum of
elements having denominator  
$(k!)^2 d_{n_1}\dots d_{n_k}$, where $n_1+\dots +n_k=2m+2$.
So we need to show that  $(k!)^2 d_{n_1}\dots d_{n_k}$ is a divisor
of $[(2m+2)!]^2 d_{2m}$ if $n_1+\dots+n_k=2m+2$.
 Since all $n_i\le 2m$, the number
$k$ is greater than 1. We assume that $n_1\le n_2\le \dots\le n_k$.
 Consider several cases. 

The case $k=2$. Then each of $n_1$, $n_2$ must be greater than 1.
By Lemma \ref{sao2}, part c),  $(k!)^2 d_{n_1}d_{n_2}= 4 d_{n_1}d_{n_2}$ is a
 divisor of
$4\times 96 d_{2m}$, which, in turn,  is a divisor of
$[(2m+2)!]^2 d_{2m}$ if $m\ge 2$.
 
The case $2m+1\ge k\ge 3$. Applying Lemma \ref{sao2}, part b), repeatedly, 
we have that 
$$d_{n_1}d_{n_2}\dots d_{n_k} \lhd 2^{k-1} 
d_{n_1+n_2+\dots+n_k-k+1}= 2^{k-1}d_{2m+3-k}.$$
If $k=3$, then $(k!)^3 d_{n_1}d_{n_2}d_{n_3}\lhd (3!)^2\times 2^2\times 
d_{2m}$, which is a divisor of $[(2m+2)!]^2 d_{2m}.$

If $k\ge 4$, then $2m+3-k\le 2m-1$. Noting that $k\le 2m+1$, one has 
$$(k!)^2 d_{n_1}\dots d_{n_k}\lhd 2^{k-1} (k!)^2 d_{2m-1}
\lhd 2^{2m} [(2m+1)!]^2 d_{2m-1}.$$
The latter is a divisor of  $[(2m+2)!]^2 d_{2m}$.

The case $k=2m+2$. Then each $n_i,i=1,2,...,2m+2$ is 1. One has
$$(k!)^2 d_{n_1}\dots d_{n_{2m+2}}= 2^{2m+2} [(2m+2)!]^2$$
is a divisor of  $[(2m+2)!]^2 d_{2m}$.
This completes the proof of Lemma \ref{denominator-psi}.
The proof of Lemma \ref{denominator-mu} is similar, even easier.

        \end{document}